\documentclass[12pt]{amsart}
 \usepackage{graphicx,mdframed}
\usepackage{amsmath,amsthm,amsfonts,amscd,amssymb,mathrsfs,hyperref}
\usepackage[all]{xypic}
\usepackage[enableskew]{youngtab} % For young tableaux
\usepackage{mathtools} %for psmallmatrix

\usepackage{color}
\usepackage[usenames,dvipsnames]{xcolor}

\newtheorem{theorem}{Theorem}[section]

\newtheorem{lemma}[theorem]{Lemma}

\newtheorem{prop}[theorem]{Proposition}

\newtheorem{computation}[theorem]{Computation}

\newtheorem{theoremst}[theorem]{Theorem$^{\star}$}

\def\A{{\mathcal A}}

\def\M{{\mathcal M}}

\newcommand{\Gr}{\operatorname{Gr}}
\newcommand{\SL}{\operatorname{SL}}

\newcommand{\GL}{\operatorname{GL}}

\newcommand{\TT}{\mathbb{T}}
\newcommand{\PP}{\mathbb{P}}

\newcommand{\CC}{\mathbb{C}}

\def\bw#1{{\textstyle \bigwedge^{\hspace{-.2em}#1}}}
\def\o{{\otimes}}

\newcommand{\VV}{V_{1}\o V_{2}\o V_{3}\o V_{4}}
\newcommand{\VVs}{V_{1}^{*}\o V_{2}^{*}\o V_{3}^{*}\o V_{4}^{*}}

\theoremstyle{definition}

\theoremstyle{remark}
\newtheorem{remark}[theorem]{Remark}

\def\red#1{{\textcolor{red}{#1}}}
\def\blue#1{{\textcolor{blue}{#1}}}
\def\green#1{{\textcolor{green}{#1}}}

\begin{document}

\author[oeding]{Luke Oeding}
\address{Department of Mathematics and Statistics, Auburn University, Auburn, AL
}
\email{oeding@auburn.edu}
\date{\today}

\title[The quadrifocal variety]{The quadrifocal variety}
\begin{abstract}
Multi-view Geometry is reviewed from an Algebraic Geometry perspective and multi-focal tensors are constructed as equivariant projections of the Grassmannian. A connection to the principal minor assignment problem is made by considering several flatlander cameras. 
The ideal of the quadrifocal variety is computed up to degree 8 (and partially in degree 9) using the representations of $\GL(3)^{\times 4}$ in the polynomial ring on the space of $3 \times 3 \times 3 \times 3$ tensors.
Further representation-theoretic analysis gives a lower bound for the number of minimal generators.
% We ask if the ideal of the quadrifocal variety is minimally generated in degree at most 9.
\end{abstract}
\maketitle
\section{Introduction and background}
\subsection{The multi-view variety}
Multi-view Geometry is a branch of Computer Vision \cite{HartleyZisserman}. An important task in Computer Vision is to efficiently reconstruct the 3-dimensional scene from the 2-dimensional projections.  Typically, one first estimates the multi-focal tensor associated to the $n$ views using correspondences arising from one object seen in multiple images.
 From the multi-focal tensor one reconstructs the camera matrices. After the camera matrices are known, one uses the point correspondences to triangulate the 3D points.

In the standard pinhole camera model, the projection of 3D world points to multiple 2D images is represented by a collection of $3\times 4$ matrices $(A_{1},\dots,A_{n})$ and the mapping
\begin{equation}\label{multi-view1}
\begin{matrix}
\PP^{3} & \to&  \PP^{2}\times \dots \times \PP^{2}
\\{}
[x] &\mapsto & ([A_{1}x],\dots,[A_{n}x]).
\end{matrix}\end{equation}
%The image of this map is the multi-view variety studied in 
For cameras in general position the multi-view mapping \eqref{multi-view1}  defines a 3 dimensional subvariety of the Cartesian product of projective spaces called the \emph{multi-view} variety.
Aholt, Sturmfels, and Thomas demonstrated rich algebraic geometry arising from this construction in \cite{AST}. They utilized a certain Hilbert scheme to describe the multi-view variety, its defining ideal, and further algebraic properties.  Their theoretical techniques  included Borel fixed monomial ideals, a universal Gr\"obner basis, degeneration to a special monomial ideal, and more.  
 This work catalyzed a new area coined ``Algebraic Vision'' by Sameer Agarwal and Rekha Thomas. The reader may wish to consult the following other examples of recent work in this field \cite{AAT,  ALST_epipolar2, ALST_epipolar1,AAT,AholtOeding,  JKSW_rigid_multiview}.

\subsection{Moduli spaces and quadrifocal tensors}Suppose  the entries of the camera matrices are not known, but are considered as parameters. By moding out by the projective rescaling in each camera plane we have a moduli space of camera matrices. The algebraic varieties of multi-focal tensors are models for these moduli spaces, and we want to know their basic algebraic properties. Bifocal tensors are just $3\times 3$ matrices of rank 2, defined by the $3\times 3$ determinant. 
Chris Aholt and the author resolved the long-standing open question of describing the ideal of trifocal tensors \cite{AholtOeding}, building on work of Alzati and Tortora \cite{Alzati-Tortora}. They also computed its algebraic degree and Hilbert polynomial using Maple, Macaulay2 \cite{M2}, and Bertini \cite{Bertini}.\looseness=-1

%We describe this construction from an algebreo-geometric point of view via equivariant projections of a Grassmannian in Section~\ref{sec:Review}.

In this paper we will be mainly concerned with the \emph{quadrifocal variety}.
One may record the correspondences induced from one 3D point seen in 4 images by the 81 special $4\times 4$ minors of  stacked camera matrix $A= (A_{1}^{\top} | A_{2}^{\top} | \cdots | A_{n}^{\top})$
 that only use one column from each of the first four blocks of $A$. These coordinates are linear in each block, yielding a $3\times 3\times 3\times 3$ \emph{quadrifocal} tensor.
  The quadrifocal variety is the Zariski closure in $\PP^{80}$ of the set of quadrifocal tensors. 

We seek a complete description of the polynomial defining equations of the quadrifocal variety. 
The main result of the present article is a first step in this direction.
\begin{theoremst}\label{thm:main} Let $I_{d}$ denote the degree $d$ piece of the ideal of the quadrifocal variety. 

$I_{d}$ is zero for  $d<3$.

$I_{3}$ is $600$-dimensional.

$I_{4}$ is $48,600$-dimensional  but contains no minimal generators.

$I_{5}$ is $1,993,977$-dimensional and contains at least $1,377$ minimal generators.

$I_{6}$ is $54,890,407$-dimensional  and contains at least $37,586$ minimal generators.

$I_{7}$ is $1,140,730,128$-dimensional  and likely contains no minimal generators.

$I_{8}$ is $18,940,147,947$-dimensional  and contains at least $162,000$ minimal generators.

$I_{9}$ is  $\geq 223,072,284,455$-dimensional  and contains at least $3,087,000$ minimal generators.
\end{theoremst}

The star refers to the fact that some of our computations were done using random points of the quadrifocal variety, so the dimensions reported are only upper bounds, but the lower bounds  hold \emph{with high probability}.
In Section~\ref{sec:Rep} we give an invariant description of all these equations and computational evidence that the equations reported here are likely \emph{not} a minimal set of generators.
%Because these are the lowest degree equations in the ideal we speculate in Conjecture~\ref{conj}, that the equations we found form a complete set of minimal generators for the ideal of the quadrifocal variety.

Until now, it was only known that a quadrifocal tensor must adhere to 51 non-linear constraints \cite{ShashuaWolf}.  
Indeed, the quadrifocal variety has codimension $51$, so there must be at least $51$ equations, but our results show that it is very far from being a complete intersection.  For instance, there are 600 cubic minimal generators.  In Section \ref{sec:Contract} we give a simple description of these equations via contractions. From the contraction description we see that the cubic equations are a consequence of the fact that every contraction of a quadrifocal tensor is a homography tensor \cite{Shashua00homographytensors}.  Additional equations are needed to take the set of tensors having that property alone  and cut it down to the quadrifocal locus.

\subsection{Multi-focal tensors in general} If $\pi$ is a partition of $4$, the entries of a multi-focal tensor of profile $\pi$ are given by the minors of $A$ that use $\pi_{i}$ columns from block $i$.
% Point or line correspondences between multiple images are encoded by these multi-focal tensors.
Multi-focal tensors record correspondences between multiple images. An \emph{epipole} is the image of one camera's focal point seen in another view. 
Correspondences between pairs of image points (in 2 views) leads to the bi-focal tensor, or \emph{fundamental matrix}. Whereas a point-point-line correspondences (in 3 views) are encoded by a \emph{trifocal tensor} \cite{PF2}, and correspondences of quadruples of image points (in 4 views) are encoded by a \emph{quadrifocal tensor} \cite{ShashuaWolf}.  See  \cite{HartleySchaffalitzky1, HartleySchaffalitzky2, ShashuaKaminski, WolfShashua, ThirthalaPollefeys} for applications of quadrifocal tensors.

Each multi-focal tensor can be determined by observing some minimal number of correspondences in multiple images.
From 7 point correspondences in 2 images one can reconstruct the bifocal tensor (fundamental matrix).  5 point correspondences in 3 or 4 images suffice to determine the trifocal and quadrifocal tensors \cite[Parts~III\&IV]{HartleyZisserman}.
For a summary of these and other minimal problems in Computer Vision, see \cite{minimal}, and in particular one can check \cite{5pt4, 5pt5, 5pt6, 4pt1, 4pt2} for recent algorithms for the relative pose problems in 3 and 4 views.
 
Since algorithms for determining the relative pose of 2 cameras exist (for instance \cite{MartinecPajdla}) one may ask what is the advantage of considering more than 2 views at a time. The first advantage is that fewer 3D points need to be identified for more than 2 views, which can be useful when many points in one view become outliers for another. In the 3 view case, the trifocal tensor can be used in the structure from motion problem when the cameras move along a straight line (see \cite{Chan2013627, VidalHartley, 4pt3}). Other advantages include greater stability, and the possibility to avoid certain unstable or critical configurations. In the 4 view case, while the quadrifocal tensor is more difficult to construct, one may avoid iterative algorithms which are not guaranteed to converge \cite[Ch.~1]{HartleyZisserman}. Moreover, unlike in the 3 view case where one view plays a special role, in the 4 view case all 4 images can play the same role, and it is expected that this symmetry can be exploited. 
For 3D to 2D projections and 4 or more views there are no $n$-focal tensors if $n>4$, hence we only consider the 2, 3 and 4 view cases.
%Image reconstruction from point correspondences: 
%Calibrated case
%N views
%P points
%6 d.o.f. per calibrated camera
%7 d.o.f for euclidean frame
%\[
%2NP \geq 3P + 6 N-7
%\]
%5 points 2 views, 4 points in 3 or more views.
%
%Uncalibrated case
%N views
%P points
%7 d.o.f. per uncalibrated camera
%7 d.o.f for euclidean frame
%\[
%2NP \geq 3P + 7 N-7
%\]
%7 points in 2 views, 5 points in 3 or 4 views. 

% Fundamental matrices are very well understood: they are $3\times 3$ matrices of rank $2$. 
%In this paper we will be mainly concerned with the quadrifocal variety, which is parameterized as follows. Collect four camera matrices in a  block matrix $A= (A_{1}^{\top} | A_{2}^{\top} | A_{3}^{\top} | A_{4}^{\top})$. The 81 coordinate functions are given by the 81 special $4\times 4$ minors of $A$ that only use one column from each block.  The quadrifocal variety is the Zariski closure of this map in $\PP^{80}$. 

%The focus of the present article is to describe the algebraic geometry of quadrifocal tensors.

\subsection{Outline}
For the reader's convenience we collect notation in Section~\ref{sec:notation}. 
We give a uniform presentation of multi-focal tensors from the viewpoint of Algebraic Geometry via equivariant projections of a Grassmannian in Section~\ref{sec:Review}, which connects our work to \cite{FaugerasMourrain,HartleySchaffalitzky1}. Section \ref{sec:Principal} addresses the case of different dimensional cameras and contains a connection between the multi-focal variety and the variety of principal minors of square matrices, connecting this work to \cite{LinSturmfels}. We discuss contractions and homography tensors in Section \ref{sec:Contract},  and state Proposition~\ref{prop:600} describing the 600 bihomogeneous cubics that are the minimal generators of the quadrifocal ideal in the lowest degree.
We use symmetry-enhanced calculations to determine the quadrifocal ideal up to degree $8$ and partially compute the ideal in degree 9 in  Section~\ref{sec:Rep}. In addition, we use representation theory and computations utilizing \texttt{SchurRings} to determine necessary minimal generators of the quadrifocal ideal. We conjecture that these necessary minimal generators also suffice. 

\section{Notation}\label{sec:notation}
Let $U$ and $V$ denote vector spaces, which we always consider the complex numbers, denoted $\mathbb{C}$, to be our ground field.
The direct sum of $U$ and $V$ is denoted $U\oplus V$, and their tensor product is denoted $U\o V$. 
The $k$-th exterior power of $V$, or the $k$-mode skew-symmetric tensors is denoted  $\bw{k}V$. Its elements are linear combinations of $k$-fold wedge products  $v_{1}\wedge v_{2}\wedge \cdots\wedge v_{n}$ with $v_{i}\in V$. We discuss alternating tensors in more depth in Section~\ref{sec:ext}.
 The vector space of symmetric $d$-mode tensors is denoted $S^{d}V$. If we choose a basis $x_{1},\ldots, x_{m}$ of $V$  we may consider $S^{d}V$ as the vector space of homogeneous degree $d$ polynomials on $x_{1},\ldots, x_{m}$. The symmetric algebra on $V$ is denoted $Sym^{\bullet}V = \bigoplus_{d\geq 0}S^{d}V$, and is isomorphic to the polynomial ring $\CC[x_{1},\ldots,x_{n}]$. The tensor algebra on $V$ is denoted $V^{\otimes } = \bigoplus_{d\geq 0} V^{\otimes d}$.

The vector space dual, denoted $V^{*}$, is the space of all linear functionals $V \to \CC$.  After bases are chosen for $U$ and $V$,  $U^{*}\o V$ may be thought of as the space of matrices representing linear mappings $U\to V$.   If $M \in U^{*}\o V$ is a matrix $\bw{k}M \in \bw{k}U^{*}\o \bw{k}V$ may be thought of as the matrix  whose entries are the $k\times k$ minors of $M$.
The general (special) linear group of all invertible (determinant 1) linear transformations of $V$ is denoted $\GL(V)$ (resp. $\SL(V)$).

We use the notation $\pi\vdash d$ to denote a partition $\pi = (p_{1},p_{2},\ldots,p_{n})$ with $\sum_{i}p_{i} = d$.
We let $S_{\pi}V$ to denote the corresponding Schur module, which we think of as an explicit representative of an irreducible submodule of the $d$-fold tensor product $V^{\o d}$.
 Given vector spaces $V_{1}, V_{2}, V_{3}, V_{4}$ and multi-partition $\pi = (\pi_{1},\pi_{2},\pi_{3},\pi_{4})$ with $\pi \vdash d$ we use the shortened notations $S_{\pi}V$ and $S_{\pi_{1}}S_{\pi_{2}}S_{\pi_{3}}S_{\pi_{4}}$ for $S_{\pi_{1}}V_{1} \o S_{\pi_{2}}V_{2}  \o S_{\pi_{3}}V_{3}  \o S_{\pi_{4}}V_{4}$.

The projective space of all lines through the origin in $V$ is denoted $\PP V$. If $v\in V$ is nonzero the line through $v$ is denoted $[v]$.
  If $X \subset \PP U $ and $Y \subset \PP V$ are algebraic varieties their Cartesian product, denoted $X\times Y$ is a subvariety of $\PP(U\oplus V)$.   The cone over the projective $X\subset \PP V$, denoted $\widehat{X}$ is an affine subvariety of $V$. 

Multifocal tensors sometimes need a mixture of indices and multi-indices. For instance, fundamental matrices are indexed by a pair of double indices: the notation $(F_{\{i,j\},\{k,l\}})$ indicates a matrix with rows indexed by the double index $\{i,j\},$ and columns indexed by the double index $\{k,l\}$, while $(T_{i,j,\{k,l\}})$ denotes a $3$-mode tensor with the first two modes indexed  by $i$ and $j$ respectively, with the third mode  indexed by the double index $\{k,l\}$. 

\section{Epipoles, fundamental matrices, trifocal and quadrifocal tensors}\label{sec:Review}

\subsection{The multi-view setup}
Multiple view geometry arises when one considers many images taken of the same scene, from (possibly) different viewpoints and is beautifully presented in \cite{HartleyZisserman}.  The following introduction is an invariant view inspired by the ideas in \cite[Ch.~17]{HartleyZisserman}. 
 Let $A_{j}$ denote $3 \times 4$ camera matrices (non-degenerate) for $1\leq j\leq n$, with row spaces equal to  $V_{j}$, a three-dimensional vector space.

Let $W$ denote a 4-dimensional vector space, whose projectivization represents the 3-dimensional ``projective world''.
For fixed camera matrices $A_{j}$, the multi view map (which also appeared in \cite{AST}) is 
\begin{equation}\label{multi-view}
\xymatrix{\PP W \ar[rr]^{\hspace{-3em}(A_{1},\dots,A_{n}) }&& \PP V_{1}\times\dots \times \PP V_{n} } 
.\end{equation}

Now we wish to treat the camera matrices as variable or as having indeterminate entries.  The map in \eqref{multi-view} is the same if we replace the matrices $A_{j}$ with scalar multiples of themselves. So, we should consider our space of parameters for cameras to be 
\[\PP (W^{*}\o V_{1}) \times \dots \times \PP (W^{*} \o V_{n}) \qquad\text{$n$-camera space}
.\]
We note here that if different camera models are taken, this paradigm may be easily altered to accommodate such changes by altering the spaces in which the cameras are modeled.

\subsection{The Grassmannian}  Faugeras and Mourrain studied multi-view geometry from the point of view of the Grassmann algebra, see \cite{FaugerasMourrain,HartleySchaffalitzky1}.  We also adapt that approach as it provides a uniform treatment and a convenient way to organize many of our computations. 

\subsubsection{Exterior products}\label{sec:ext}
Recall if $V$ is a vector space with basis $u_{1},\dots,u_{n}$, we construct the exterior powers of $U$, denoted $\bw{k}U$ by considering the alternating (or wedge) product $\wedge$ and forming the vector space of $k$-vectors (length $k$ wedge products) with basis consisting of pure $k$-vectors $\left \{u_{i_{1}},\dots,u_{i_{k}}\mid 1\leq i_{1}<i_{2}<\dots<i_{k}\leq n \right\}$. 
Thus $\bw{k} U$ has dimension $\binom{n}{k}$ if $0\leq k\leq n$ and 0 otherwise.  It is straightforward to see that  $v_{1}\wedge\dots\wedge v_{k} \neq 0$  if and only if the vectors $\{v_{1},\ldots,v_{k}\}$ are linearly independent in $U$. Consider two non-zero $k$-vectors $v_{1}\wedge\dots \wedge v_{k}$ and $w_{1}\wedge\dots \wedge w_{k}$ and the underlying vector spaces $E:=\textrm{span}\{v_{1},\dots, v_{k}\}$ and $F:=\textrm{span}\{w_{1},\dots, w_{k}\}$. It is straightforward to check that 
\[
v_{1}\wedge\dots\wedge v_{k} = \lambda ( w_{1}\wedge\dots \wedge w_{k}) \quad \textrm{for some} \quad \lambda \neq0 \quad \textrm{if and only if} \quad \bw{k}E = \bw{k}F,
\]
and equality holds when $\lambda$ is the determinant of the change of basis between $E$ and $F$.

We denote by $\PP \bw{k} U$ the projective space consisting of lines through the origin in $\bw{k} U$, which we may consider as the set of classes $[\omega] = \left\{\lambda \omega \mid \lambda \in\CC \setminus \{ 0\}, \omega \in \bw{k}U \setminus \{ 0\} \right \}$.

This motivates the definition of the Grassmannian (in its minimal embedding). Let $\Gr(k,U)$ denote the set of $k$-planes in $U$.   The rational map:
\[\begin{matrix}
\Gr (k,U) &\to& \PP \bw{k} U\\
 M & \mapsto & \bw{k}M
 \end{matrix}
\]
is an embedding (in fact, the embedding is a minimal rational embedding). The usual Pl\"ucker embedding is a slight variant of this construction, which we will review next in the context of multiple view geometry.

\subsubsection{From multiple views to the Grassmannian}
It is natural to consider the following  $4\times 3n$ blocked matrix, which will present a convenient way to keep track of external constraints on the multi-view setup
\[
M = \left(\begin{array}{c|c|c|c} A_{1}^{\top} & A_{2}^{\top}&\dots & A_{n}^{\top}\end{array}\right) \in (W^{*}\o V_{1})\oplus \dots\oplus (W^{*}\o V_{n}) = W^{*}\o (\textstyle \bigoplus_{j}V_{j})
.\]
The non-degeneracy condition is that each matrix $A_{j}$ in an $n$-tuple $([A_{1}^{\top}],[A_{2}^{\top}],\dots,[A_{n}^{\top}])$ must have full rank, which occurs in an open set; so the row space of $M$ parameterizes the (Grassmannian variety) 4-dimensional subspaces of a $3n$-dimensional space.

The maximal minors of $M$ give coordinates (the Pl\"ucker coordinates) on the Grassmannian $\Gr(4,3n)$. 
These minors are also known as \emph{multilinear coordinates}.

In invariant language this parameterization is the following map
\begin{equation}\label{plucker}
\begin{array}{rcl}\varphi\colon 
 W^{*}\o (\bigoplus_{j}V_{j}) &\longrightarrow  &\bw{4}W^{*}\o \bw{4}(\bigoplus_{j}V_{j}) \cong \bw{4}(\bigoplus_{j}V_{j})
 \\ M &\longmapsto & \bw{4}M 
\end{array}
.\end{equation}
The image of $\varphi$ (the Zariski closure of the image of an open set) is isomorphic to the cone over the Grassmannian of 4-dimensional subspaces of $\bigoplus_{j}V_{j}$, which is the row space of $M$. In other words
\[\textstyle
Im(\varphi) =
\widehat{\Gr}(4, \bigoplus_{j}V_{j}) \subset  \bw{4}(\bigoplus_{j}V_{j}).
\]
Notice that because $W$ is 4-dimensional, $\bw{4}W$ is one-dimensional, and thus passing to the maximal minors of the concatenated matrix $M$ removes the dependency on the world points represented by $\PP W$. 

It is well known that the dimension of the Grassmannian $\Gr(r,\CC^{N})$ is $r(N-r)$. In our example $\dim(\Gr(4,3n)) = 4(3n-4)$.

If we restrict to camera space, we have to consider the image of the map up to the $n$-dimensional torus action which records the projective ambiguity in each of the $n$ cameras. We can restrict the target of the map to the appropriate \emph{GIT quotient}:
\[\textstyle
\PP (W^{*}\o V_{1}) \times \dots \times \PP (W^{*} \o V_{n}) \to \widehat{Gr}(4, \bigoplus_{j}V_{j}) /\!/ (\CC^{*})^{n} \subset  \bw{4}(\bigoplus_{j}V_{j})/\!/ (\CC^{*})^{n}.
\]
From this we obtain the dimension of the GIT quotient  (see \cite[\S~17.5]{HartleyZisserman} and \cite[\S~6]{AST}) 
\[\dim \widehat{Gr}(4, \bigoplus_{j}V_{j}) /\!/ (\CC^{*})^{n} \quad=\quad 4(3n-4) +1 - n \quad =\quad 11n-15.\]

\begin{remark}Here is a classical formula for the degree of the Grassmann manifold, (see \cite{MukaiGrassmannian},  \cite[Ex.~14.7.11]{FultonINT} or \cite[\S~10.1.2, eq. (10.6)]{DolgachevAG}) 
\[\deg\Gr(r,n)  \quad = \quad 
\frac{
1!2!\cdots r! \dim(\Gr(r,n))!
}
{
(n-r)!(n-r+1)! \cdots n!
}
%(r(n - r))! \prod\limits_{1\leq i \leq r \leq j \leq n} (j - i)^{-1}
.\]
For example
\[
\deg(\Gr(4,6)) = 14, \qquad
\deg(\Gr(4,9)) = 1662804, \qquad
\deg(\Gr(4,12)) = 1489877926680.
\]
One hope is that a better understanding of the GIT quotient of the Grassmannian might allow us to find the degree of the multi-focal tensor varieties.
\end{remark}

\subsection{Symmetry}
For $j \in \{1,2,3,4\}$ let $A_{j}$ be a $3 \times 4$ (non-degenerate) camera matrix (an element of $W^{*}\o V_{j}$), with blocking $A_{j} = (B_{j}| {\bf x}_{j})$.  
On each matrix $A_{j}$ we have an action of $\GL(V_{j}) \cong \GL(3)$ acting by change of coordinates in the camera plane. The action is
\begin{multline*}
\left(\GL(V_{1})\times \GL(V_{2})\times \GL(V_{3})\times \GL(V_{4}) \right)
\times 
\left( (W^{*}\o V_{1}) \oplus (W^{*}\o V_{2}) \oplus (W^{*}\o V_{3}) \oplus (W^{*}\o V_{4})  \right)\\
\to \left( (W^{*}\o V_{1}) \oplus (W^{*}\o V_{2}) \oplus (W^{*}\o V_{3}) \oplus (W^{*}\o V_{4})  \right)\\
(g_{1},g_{2},g_{3},g_{4}),(A_{1},A_{2},A_{3},A_{4}) \mapsto (g_{1}A_{1},g_{2}A_{2},g_{3}A_{3},g_{4}A_{4}),
\end{multline*}
where we take the action of each $g_{j}$ to be a change of basis in the row space of $A_{j}$.
Because the matrices $A_{j}$ are assumed to be full rank, we can, without loss of generality, act by an element of $GL(3)^{\times 4}$ and  assume that $B_{j}= Id_{3}$ and move the 4-tuple $(A_{1}^{\top} | A_{2}^{\top} | A_{3}^{\top} | A_{4}^{\top})$ to 
\begin{equation}\label{Amatrix}
A \cong \left( \begin{array}{c|c|c|c}
 \begin{smallmatrix}
1 & 0 & 0 \\
0 & 1 & 0 \\
0 & 0 & 1 \\
x_{1,1}& x_{1,2} & x_{1,3}
\end{smallmatrix}
&
 \begin{smallmatrix}
1 & 0 & 0 \\
0 & 1 & 0 \\
0 & 0 & 1 \\
x_{2,1}& x_{2,2} & x_{2,3}
\end{smallmatrix}
&
 \begin{smallmatrix}
1 & 0 & 0 \\
0 & 1 & 0 \\
0 & 0 & 1 \\
x_{3,1}& x_{3,2} & x_{3,3}
\end{smallmatrix}
&
 \begin{smallmatrix}
1 & 0 & 0 \\
0 & 1 & 0 \\
0 & 0 & 1 \\
x_{4,1}& x_{4,2} & x_{4,3}
\end{smallmatrix}
\end{array}\right)
.\end{equation}

\begin{remark}
There is an action of $\GL(W)$ acting on simultaneously on all of the column spaces of $A_{j}$ (which are all equal to $W$). While this action doesn't turn out to be useful (it is trivialized in the tensorial map), the action of $\mathfrak{S}_{4}$ permuting the matrices $A_{j}$, also permutes the indices in the image of the tensorial map, and preserves the set of quadrifocal tensors.
\end{remark}

\subsection{Multilinearity spaces}
The $4\times 4$ minors of $M$ come in several classes, which have invariant descriptions.
By construction, 
$ \bw{4}(\bigoplus_{j}V_{j})$ is a vector space with a natural $\GL(3n)$-action, but we can further view this as a $\GL( \bw{4}(\bigoplus_{j}V_{j}))$-action (each $\GL(V_{j})$ acting by invertible linear change of coordinates in $V_{j}$), and there is a natural inclusion of $G:=\prod_{j}\GL(V_{j})$, which may be thought of as block diagonal (with the proper choice of basis) inside $\GL(3n)$.  Moreover, we may view $G$ as a product of a torus $\TT^{n}:=\prod_{j}T_{j}\cong (\CC^{*})^{n}$ ($T_{j}\cong\CC^{*}$ acting by scaling block $j$ of $M$) and $\prod_{j}\SL(V_{j})$. In summary
\[G := \prod_{j}T_{j} \times \prod_{j}\SL(V_{j}) \subset \GL(\bw{4}\bigoplus V_{j}).\]

Thus, we may consider $\Gr(4,\bigoplus V_{j})$ as a $G$-variety. On the other hand, on the GIT quotient, the torus action is trivialized, so we will consider $G':=\prod_{j}\SL(V_{j})$ acting on $\Gr(4, \bigoplus V_{j})\subset \bw{4}\bigoplus V_{j}$ and on the GIT quotient $ \widehat{Gr}(4, \bigoplus_{j}V_{j}) /\!/ (\CC^{*})^{n} \subset  \bw{4}(\bigoplus_{j}V_{j})/\!/ (\CC^{*})^{n}$.  The effect of trivializing the torus action is that we may identify every $V_{j}$ with its dual, and this induces an identification of every irreducible representation $S_{\pi}V_{j}$ with its dual $S_{\pi}V_{j}$.

Now $G$ acts on $\bw{4}(\bigoplus_{j}V_{j})$ and it has a decomposition into irreducible $G$-modules:
\[\begin{matrix}
\bw{4}\left(\bigoplus_{j}V_{j}\right) &= &
\left(  \bigoplus_{i\neq j} \bw{3}V_{i}\o V_{j} \right) \oplus
\left(  \bigoplus_{i < j} \bw{2}V_{i}\o \bw{2} V_{j} \right) \\
&&\oplus
\left(  \bigoplus_{i\not\in \{j,k\}, j<k} \bw{2}V_{i}\o V_{j} \o V_{k} \right) \oplus
\left(  \bigoplus_{i<j<k<l}V_{i}\o V_{j} \o V_{k} \o V_{l} \right) .\end{matrix}
\]
%where the summations are over (respectively) pairs, triples and quadruples of  distinct indices.
The 4 non-isomorphic module classes and their descriptions are listed in Table~\ref{table1}.
\begin{table}%[h!]
\[
\begin{array}{|rrrl|c|l|}
\hline
&&&\hspace{-5em}\text{Columns used}  &\text{Invariant description}  & \text{Space}  \\
\hline
i_{1}&i_{2}&i_{3}&j &         \bw{3}V_{i}\o V_{j} \;\;\cong\;\; V_{j} &\text{\emph{epipole space}} \\
i_{1}&i_{2}&j_{1}&j_{2}&  \bw{2}V_{i}\o \bw{2} V_{j} \;\;\cong\;\; V_{i}^{*}\o V_{j}^{*} & \text{\emph{fundamental matrix space}}\\
i_{1}&i_{2}&j&k&              \bw{2}V_{i}\o V_{j}\o V_{k} \;\;\cong\;\; V_{i}^{*}\o V_{j}\o V_{k}  &\text{\emph{trifocal space}}\\
i&j&k&l&                           V_{i}\o V_{j}\o V_{k}\o V_{l}   & \text{\emph{quadrifocal space}}\\
\hline
\end{array}
\]
\caption[Table]{
\emph{Classes of multi-focal tensors.} Assume $i,j,k,l$ are distinct,  $i_{s}$ are in block $i$, and  $j_{s}$ are in block $j$.
}\label{table1}
\end{table}

Now we will consider the projection of the cone over the Grassmannian $\widehat{Gr}(4,\bigoplus_{j}V_{j})$ to each type of multi-linearity space.  The images of the projections are  respectively the single view, the epipolar variety, the trifocal variety, and the quadrifocal variety.  The fiber of the projection over a general point is the product of the ignored camera planes and the torus acting on the utilized camera planes. So it is interesting to consider the minimal number of cameras in each case.  Moreover, because the projections on the level of vector spaces are equivariant, the images are automatically invariant (with respect to the appropriate group).

Suppose now that there are at most $4$ cameras. For each $i \in \{1,2,3,4\}$ let the set $\{e^{i}_{1},e^{i}_{2},e^{i}_{3}\}$ denote a basis of $V_{i}$, which also provides an ordered basis on  $\CC^{12}\cong V_{1}\oplus V_{2}\oplus V_{3}\oplus V_{4}$. 
In the remainder of this section we consider the parameterization of 
epipoles Sec.~\ref{epi},
fundamental matrices Sec.~\ref{fun},
trifocal tensors Sec.~\ref{tri},
and quadrifocal tensors Sec.~\ref{quad}, all up to the $G$-action using the matrix $A$ in \eqref{Amatrix}. These results are summarized in Table~\ref{table2}.

\begin{table}%[h!]
\[
\begin{array}{|c|c|}
\hline & \\[-10pt]
S = \begin{psmallmatrix}
( x_{1,1}-x_{2,1}) &
(-1)( x_{1,2}-x_{2,2})&
( x_{1,3}-x_{2,3})
\end{psmallmatrix}^{\top}
&\text{\emph{$(1,2)$-epipole}} \\
&\\[-10pt]
\hline &\\[-10pt]
F = \begin{psmallmatrix}0&
      x_{1,3}-x_{2,3}&
      -x_{1,2}+x_{2,2}\\
      -x_{1,3}+x_{2,3}&
      0&
      x_{1,1}-x_{2,1}\\
      x_{1,2}-x_{2,2}&
      -x_{1,1}+x_{2,1}&
      0\\
\end{psmallmatrix}
      &\text{\emph{$(1,2)$-fundamental matrix}}\\
&\\[-10pt]      \hline&\\[-10pt]
\begin{array}{ll}
\left. \begin{smallmatrix}T_{i,i,\{k,l\}}  &=&  -x_{1,i}+x_{2,i} \\
T_{i,k,\{i,l\}}  &=&  -x_{1,i}+x_{3,i} \\
T_{k,i,\{i,l\}}  &=&  x_{2,i}-x_{3,i} ,
\end{smallmatrix} \right\}& \begin{array}{l}k, l \text{ distinct, } \\i  \text{ distinct from } k \,\&\, l,\end{array} \\
\left.\begin{smallmatrix} T_{i,j,\{k,l\}} &=& 0 \end{smallmatrix}\hspace{3.2em} \right\}& \text{ else}.
\end{array}
&\text{\emph{$(1,2,3^{*})$-trifocal tensor}}\\
&\\[-10pt]
\hline
&\\[-10pt]
\begin{array}{ll}
\left. \begin{smallmatrix}
Q_{i,i,k,l} &=& (-1)^{i} (x_{1,i}-x_{2,i})  &=&  -Q_{i,i,l,k},\\
Q_{i,k,i,l}  &=& (-1)^{i} (x_{1,i}-x_{3,i})  &=&  -Q_{i,l,i,k}, \\
Q_{i,k,l,i}  &=&  (-1)^{i} (x_{1,i}-x_{4,i})  &=&  -Q_{i,l,k,i}, \\
Q_{k,i,i,l}  &=&  (-1)^{i} (x_{2,i}-x_{3,i})  &=& -Q_{l,i,i,k},\\
Q_{k,i,l,i}  &=&  (-1)^{i} (x_{2,i}-x_{4,i})  &=& -Q_{l,i,k,i},\\
Q_{k,l,i,i}  &=&  (-1)^{i} (x_{3,i}-x_{4,i})  &=&  -Q_{l,k,i,i}
\end{smallmatrix} \right\}& \begin{array}{l}k< l \text{,  } \\ i  \text{ distinct from } k \,\&\, l,  \end{array} \\
\left.\begin{smallmatrix}Q_{i,j,k,l} &=& 0 \end{smallmatrix}\hspace{9em} \right\} & \text{ else}.
\end{array}
 & \text{\emph{$(1,2,3,4)$-quadrifocal tensor}}\\
\hline
\end{array}
\]
\caption[Table]{\emph{Parametrizations of multi-focal tensors up to symmetry.}}\label{table2}
\end{table}

\subsection{Epipoles}\label{epi}
For a pair of cameras, we consider the projection of $\Gr(4,6) = \Gr(4,V_{1}\oplus V_{2})\subset \PP \bw{4}(V_{1}\oplus V_{2}) = \PP^{14}$ to epipolar space $\PP (V_{1}\o \bw{3}V_{2}) = \PP^{2}$. The target space is (naturally isomorphic to) the projective plane, and the map subjects onto $\PP^{2}$. 
 The image is naturally $\GL(V_{1})\times \GL(V_{2})$-invariant, the action of $\GL(V_{1})$ being trivial and the 2-dimensional torus acts by a weight of $(3,1)$. 
 
% We also may consider the maps on the level of the GIT quotients: $\widehat{\Gr(4,6)}/\!/ (\CC^{*})^{2}$  mapping to $( V_{1}\o \bw{3}V_{2})/\!/ (\CC^{*})^{2} \cong \PP V_{2}$. 

The image of the projection is the space of epipoles in view 1 imposed by view 2.  Thus the epipoles may be recovered from the multi-view setup via a projection from the Grassmannian.
To get an expression of an epipole in coordinates  consider just two cameras and the matrix
\[
A \cong \left( \begin{array}{c|c}
 \begin{smallmatrix}
1 & 0 & 0 \\
0 & 1 & 0 \\
0 & 0 & 1 \\
x_{1,1}& x_{1,2} & x_{1,3}
\end{smallmatrix}
&
 \begin{smallmatrix}
1 & 0 & 0 \\
0 & 1 & 0 \\
0 & 0 & 1 \\
x_{2,1}& x_{2,2} & x_{2,3}
\end{smallmatrix}
\end{array}\right)
.\]
An element $S$ of single camera space has the following form  in the Pl\"ucker coordinates associated to 
the determinants of the matrices constructed from  one column from the first block of $A$ the three columns of the second block:
\[S = \sum_{1\leq p \leq 3} S_{p,\{1,2,3\}}e^{1}_{p}\wedge e^{2}_{1}\wedge e^{2}_{2} \wedge e^{2}_{3}
\quad \in \quad    \bw{3}V_{2}\o V_{1} \cong V_{1}
.\] 
Applying this to $A$ we get the coordinates of the epipole
\[\begin{pmatrix}
S_{1,\{1,2,3\}}(A) \\
S_{2,\{1,2,3\}}(A) \\
S_{3,\{1,2,3\}}(A) 
\end{pmatrix}
=
\begin{pmatrix}
( x_{1,1}-x_{2,1})\\
(-1)( x_{1,2}-x_{2,2})\\
( x_{1,3}-x_{2,3})
\end{pmatrix}
.\]

\subsection{Fundamental matrices}\label{fun}
Again for a pair of cameras, we may consider the projection of $\Gr(4,6) = \Gr(4,V_{1}\oplus V_{2})\subset \PP \bw{4}(V_{1}\oplus V_{2}) = \PP^{14}$ to fundamental matrix space $\PP (V_{1}^{*}\o V_{2}^{*}) = \PP^{8}$. The target space may be interpreted as the projectivization of a space of $3 \times 3$ matrices, and caries the natural action of $\GL(V_{1})\times \GL(V_{2})$. We might call the image of the projection the \emph{variety of fundamental matrices}, which is also naturally $\GL(V_{1})\times \GL(V_{2})$-invariant. Because the vector spaces $V_{1}$ and $V_{2}$ play symmetric roles, this image is also naturally $\mathfrak{S}_{2}$ invariant.
It is well known that the matrices in the image of the projection have a one-dimensional kernel and the image variety is just the (degree 3) determinantal hypersurface.

Note that $\widehat{Gr}(4,6)$ is 9-dimensional, the GIT quotient $\widehat{\Gr(4,6)}/\!/ (\CC^{*})^{2}$ is 7-dimensional, and thus birational to the projective variety of $3 \times 3$ matrices of rank $\leq 2$, the variety of fundamental matrices.  Also note that the 2-dimensional torus acts by a weight of $(2,2)$.

An element $F$ of fundamental matrix space has the following form  in Pl\"ucker coordinates:
\[F = \sum_{1\leq i,j,k,l \leq 3,\; i<j, \;k<l} F_{\{i,j\},\{k,l\}}e^{1}_{i}\wedge e^{1}_{j}\wedge e^{2}_{k} \wedge e^{2}_{l}
\quad \in \quad \bw{2}V_{i}\o \bw{2} V_{j} \cong V_{i}^{*}\o V_{j}^{*}
.\]
For a fixed pair of distinct indices $i,j$, a fundamental matrix is gotten by applying $F$ to a blocked $4\times 12$ camera matrix $A$. In particular, $F(A)$ may be thought of as a vector whose coordinates are determinants of the $4\times 4$ submatrices of $A$ obtained by taking two columns from each from blocks $i$ and $j$ of $A$.

Now apply $F$ to $A$ in \eqref{Amatrix}. The (1-2) fundamental matrix associated to $A$ is described by
coordinates $ F_{\{i,j\},\{k,l\}}$, and we can represent this as the matrix
\[F(A) = 
\begin{pmatrix}0&
      x_{1,3}-x_{2,3}&
      -x_{1,2}+x_{2,2}\\
      -x_{1,3}+x_{2,3}&
      0&
      x_{1,1}-x_{2,1}\\
      x_{1,2}-x_{2,2}&
      -x_{1,1}+x_{2,1}&
      0\\
      \end{pmatrix},\] where the $(p,q)$ entry is $ F_{\{i,j\},\{k,l\}}$ such that $\{p,i,j\} = \{q,k,l\} =\{1,2,3\}$. 

\subsection{The trifocal variety} \label{tri}
For a triple of cameras, we consider the projection of $\Gr(4,9) = \Gr(4,V_{1}\oplus V_{2} \oplus V_{3})\subset \PP \bw{4}(V_{1}\oplus V_{2} \oplus V_{3}) = \PP^{\binom{9}{4}-1}$ to trifocal space $\PP (V_{1}\o V_{2}\o V_{3}^{*}) = \PP^{26}$. The target space may be interpreted as the projectivization of a space of $3 \times 3 \times 3$ tensors, and caries the natural action of $\GL(V_{1})\times \GL(V_{2})\times \GL(V_{3})$. The image of the projection is called the \emph{trifocal variety}, which is also naturally $\GL(V_{1})\times \GL(V_{2})\times \GL(V_{3})$-invariant.
Because the vector spaces $V_{1}$ and $V_{2}$ play symmetric roles, this image is also naturally $\mathfrak{S}_{2}$ invariant. Non-trivial permutations involving $V_{3}$ will not preserve this trifocal variety but produce an isomorphic copy it.

The trifocal variety is 18-dimensional (\cite[p.~368]{HartleyZisserman}). The GIT quotient $\widehat{\Gr(4,9)}/\!/ (\CC^{*})^{3}$ is 18-dimensional and thus is birational to the trifocal variety.  In this case, the 3-dimensional torus acts by a weight of $(1,1,2)$.

For the trifocal tensor, consider 3 cameras in special position producing 
\[
A \cong \left( \begin{array}{c|c|c}
 \begin{smallmatrix}
1 & 0 & 0 \\
0 & 1 & 0 \\
0 & 0 & 1 \\
x_{1,1}& x_{1,2} & x_{1,3}
\end{smallmatrix}
&
 \begin{smallmatrix}
1 & 0 & 0 \\
0 & 1 & 0 \\
0 & 0 & 1 \\
x_{2,1}& x_{2,2} & x_{2,3}
\end{smallmatrix}
&
 \begin{smallmatrix}
1 & 0 & 0 \\
0 & 1 & 0 \\
0 & 0 & 1 \\
x_{3,1}& x_{3,2} & x_{3,3}
\end{smallmatrix}
\end{array}\right)
.\]
An element $T$ of trifocal space has the following form  in Pl\"ucker coordinates:
\[T =  \sum_{1\leq i,j,k,l \leq 3,\; k<l}
T_{i,j,\{k,l\}} e^{1}_{i}\wedge e^{2}_{j}\wedge e^{3}_{k} \wedge e^{3}_{l}
\quad \in \quad 
V_{1}\o V_{2}\o   \bw{2}V_{3} \cong V_{1}\o V_{2}\o V_{3}^{*}
.\] 
A trifocal tensor is gotten by applying $T$ to  $A$. In particular, $T(A)$ may be thought of as a vector whose coordinates are determinants of the $4\times 4$ submatrices of $A$ obtained by taking one column from each of the first two blocks of $A$ and two columns of from the third block of $A$.

\begin{lemma}
The locus of trifocal tensors is the $\GL(3)\times\GL(3)\times \GL(3)$-orbit of the 9-dimensional linear space whose coordinates are given by $T_{i,j,\{k,l\}}$ satisfying the following conditions:
\[\begin{array}{ll}
\left. \begin{array}{l}T_{i,i,\{k,l\}}(A) = -x_{1,i}+x_{2,i} \\
T_{i,k,\{i,l\}}(A) = -x_{1,i}+x_{3,i} \\
T_{k,i,\{i,l\}}(A) = x_{2,i}-x_{3,i} ,
\end{array} \right\}& k \text{ and } l \text{ distinct, and } i  \text{ distinct from } k,l, \\
\begin{array}{l} T_{i,j,\{k,l\}}(A)=0 \end{array} & \text{else}.
\end{array}
\]
\end{lemma}

\subsection{The quadrifocal variety}\label{quad}
For a quadruple of cameras, we consider the projection of $\Gr(4,12) = \Gr(4,V_{1}\oplus V_{2} \oplus V_{3}\oplus V_{4})\subset \PP \bw{4}(V_{1}\oplus V_{2} \oplus V_{3}\oplus V_{4}) = \PP^{\binom{12}{4}-1}$ to quadrifocal space $\PP (V_{1}\o V_{2}\o V_{3}\o V_{4}) = \PP^{80}$. The target space may be interpreted as the projectivization of a space of $3 \times 3 \times 3 \times 3$ tensors, and caries the natural action of $\GL(V_{1})\times \GL(V_{2})\times \GL(V_{3})\times \GL(V_{4})$. The image of the projection is called the \emph{quadrifocal variety}, which is also naturally $\GL(V_{1})\times \GL(V_{2})\times \GL(V_{3})\times \GL(V_{4})$-invariant.
Because the vector spaces $V_{1},V_{2}, V_{3}$ and $V_{4}$ play symmetric roles, this image is also naturally $\mathfrak{S}_{4}$ invariant. 

The quadrifocal variety is 29-dimensional (\cite[p.~423]{HartleyZisserman}), and the GIT quotient $\widehat{\Gr(4,12)}/\!/ (\CC^{*})^{4}$ is 29-dimensional, and thus birational to the quadrifocal variety.  In this case, the 4-dimensional torus acts by a weight of $(1,1,1,1)$.

An element $Q$ of quadrifocal space has the following form in Pl\"ucker coordinates:
\[Q = \sum_{1\leq i,j,k,l \leq 3}Q_{i,j,k,l}e^{1}_{i}\wedge e^{2}_{j}\wedge e^{3}_{k} \wedge e^{4}_{l}
\quad \in \quad  V_{i}\o V_{j}\o V_{k}\o V_{l}   
.\]
A quadrifocal tensor is gotten by applying $Q$ to a blocked $4\times 12$ camera matrix $A$. In particular, $Q(A)$ may be thought of as a vector whose coordinates are determinants of the $4\times 4$ submatrices of $A$ obtained by taking one column from each of the 4 blocks of $A$.

%\begin{lemma} Consider $A$ in \eqref{Amatrix}. In Pl\"ucker coordinates,  we have the Quadrifocal Space associated to $A$ described by
%\[
%Q(A)_{1,2,3,4} = \sum_{
%\begin{matrix}1\leq i_{1},i_{2},i_{3},i_{4} \leq 3
%\\
%i_{p}=i_{q}=r
%\end{matrix}
%} sgn(r; i_{1},i_{2},i_{3},i_{4})
% (x_{p,r}-x_{q,r}) e^{1}_{i_{1}} \wedge e^{2}_{i_{2}}\wedge e^{3}_{i_{3}}\wedge e^{4}_{i_{4}}
%\]
%
% \[Q(A)_{1,2,3,4} = T(A)_{1,2,\ul{3,4}} + T(A)_{1,3,\ul{2,4}} + T(A)_{1,4,\ul{2,3}} % + T(A)_{2,3,\ul{1,4}} + T(A)_{2,4,\ul{1,3}} + T(A)_{3,4,\ul{1,2}}.
% \]
% In particular, we have an expression of Quadrifocal Space as a $\GL(3)\times \GL(3)\times \GL(3)\times\GL(3)$-orbit of the 12-dimensional vector space parameterized by $\{x_{1,i},x_{2,j}, x_{3,k},x_{4,l}\}$.
%\end{lemma}
%
%
%
%\subsubsection{A normal form for the quadrifocal tensor}
For the quadrifocal tensor, consider 4 cameras in special position.
\[
A \cong \left( \begin{array}{c|c|c|c}
 \begin{smallmatrix}
1 & 0 & 0 \\
0 & 1 & 0 \\
0 & 0 & 1 \\
x_{1,1}& x_{1,2} & x_{1,3}
\end{smallmatrix}
&
 \begin{smallmatrix}
1 & 0 & 0 \\
0 & 1 & 0 \\
0 & 0 & 1 \\
x_{2,1}& x_{2,2} & x_{2,3}
\end{smallmatrix}
&
 \begin{smallmatrix}
1 & 0 & 0 \\
0 & 1 & 0 \\
0 & 0 & 1 \\
x_{3,1}& x_{3,2} & x_{3,3}
\end{smallmatrix}
&
 \begin{smallmatrix}
1 & 0 & 0 \\
0 & 1 & 0 \\
0 & 0 & 1 \\
x_{4,1}& x_{4,2} & x_{4,3}
\end{smallmatrix}
\end{array}\right)
.\]

\begin{lemma}\label{lem:quadForm}
The locus of quadrifocal tensors is the $\GL(3)\times\GL(3)\times \GL(3)\times \GL(3)$-orbit of the 12-dimensional linear space of tensors whose coordinates are given by $Q_{i,j,k,l}$ satisfying the following conditions:
\[\begin{array}{ll}
\left. \begin{array}{l}
Q_{i,i,k,l}(A) =(-1)^{i} (x_{1,i}-x_{2,i}) = -Q_{i,i,l,k}(A),\\
Q_{i,k,i,l}(A) =(-1)^{i} (x_{1,i}-x_{3,i}) = -Q_{i,l,i,k}(A), \\
Q_{i,k,l,i}(A) = (-1)^{i} (x_{1,i}-x_{4,i}) = -Q_{i,l,k,i}(A), \\
Q_{k,i,i,l}(A) = (-1)^{i} (x_{2,i}-x_{3,i}) =-Q_{l,i,i,k}(A),\\
Q_{k,i,l,i}(A) = (-1)^{i} (x_{2,i}-x_{4,i}) =-Q_{l,i,k,i}(A),\\
Q_{k,l,i,i}(A) = (-1)^{i} (x_{3,i}-x_{4,i}) = -Q_{l,k,i,i}(A)
\end{array} \right\}& k< l \text{  and } i  \text{ distinct from } k,l, \\
\begin{array}{l} Q_{i,j,k,l}(A)=0 \end{array} & \text{else}.
\end{array}
\]
\end{lemma}

\section{Arbitrary dimensional cameras: a connection to principal minors}\label{sec:Principal}
It is curious to study the case when $W$ and  $V_{i}$ arbitrary dimensions (see \cite{BertoliniBesanaTurrini}).  For instance, when each $V_{i}$ has dimension two, we might imagine the cameras to be flatlanders' cameras.  When the $V_{i}$ have dimension more than 3, we might consider an image plane that records higher dimensional data such as temperature or color. 

\subsection{Principal minors}
In the case $V_{i}\cong\CC^{2}$ we find it very interesting that seemingly unrelated work regarding principal minors fits into this framework \cite{HoltzSturmfels,LinSturmfels,BorodinRains,Oeding_principal, oeding_thesis,GriffTsat1,GriffTsat2}.

Suppose $W$ has dimension $m$ and that we have $m$ copies of $V_{i}$, each with dimension 2.
Consider again the matrix $M = \left(\begin{array}{c|c|c|c} A_{1}^{\top} & A_{2}^{\top}&\dots & A_{m}^{\top}\end{array}\right)$, now as a $m\times 2m$ matrix consisting of $m$ blocks of size $m\times 2$.  By a left action of $\GL(W)$ we may assume that each block of $M$ is of the form $A_{i}^{\top} = \begin{pmatrix}e_{i}& b_{i}\end{pmatrix}$, where $e_{i}$ is the $i$-th standard basis vector of $W$ and $b_{i}$ is arbitrary. Let $B$ denote the $m\times m$ matrix with columns $b_{i}$.

It is straightforward to see that the maximal minors of $M$ correspond to the minors of $B$ and  moreover that the minors of $M$ that use precisely one column from each block of $M$ correspond to the principal minors of $B$.

In turn, this identification naturally gives the space $V_{1}\o\dots\o V_{m}$ the interpretation as the space of all $2^{m}$ principal minors of an $m\times m$ matrix (the $0\times 0$ minor may be assumed to be $1$ in this construction). 

The case $m=4$ was studied by  \cite{BorodinRains, LinSturmfels}, who discovered the minimal defining equations of the ideal of relations among principal minors of a generic $4\times 4$ matrix, and found that the algebraic variety coincides with the main component of the singular locus of the $2\times 2\times 2\times 2$ hyperdeterminant. 

\begin{remark}
Weyman and Zelevinski \cite{WeyZel_Sing} analyzed the singular locus of the hyperdeterminant. In the $2\times2\times2\times 2$ case, there are 8 components of the singular locus:
\[
\nabla_{\text{cusp}} \cup \nabla_{\text{node}}(\emptyset) \cup \nabla_{\text{node}}(\{i,j\}) \quad 1\leq i< j \leq 4.
\]
Moreover, according to Holweck, \cite{Holweck_SingGrass} the ``main'' component $\nabla_{\text{node}}(\emptyset)$ also has the interpretation as the projective dual of the secant line variety:
\[
\nabla_{\text{node}}(\emptyset) = \sigma_{2}(\PP^{1}\times \PP^{1}\times \PP^{1}\times \PP^{1})^{\vee} \subset \PP (\VV)^{*}
.\]
The variety of principal minors of $4\times 4$ matrices is the projection 
\[
\Gr(4,V_{1}\oplus V_{2} \oplus V_{3} \oplus V_{4}) \to \PP(\VV)
.\]
The dual variety of $\Gr(4,8)$ is \emph{not} one of the 3 cases which is normal \cite{Holweck_SingGrass}, however it does have many other nice properties, such as a finitely generated invariant ring, finitely many orbits in the null cone, and a set of normal forms that depend on 8 parameters, \cite{CDZG}.   
It would be nice to see how to exploit these coincidences. 
\end{remark}

In the case that the matrix $B$ is assumed to be symmetric, \cite{HoltzSturmfels} re-invigorated a classical study going back to Cayley of relations amongst principal minors of symmetric matrices. (See \cite{LinSturmfels} and \cite{Oeding_principal} for a historical discussion.)
In particular, the ideal of relations of the principal minors of a generic symmetric $3 \times 3$ matrix is generated by a single equation known as Cayley's $2\times 2\times 2$ hyperdeterminant, denoted $h$ for this discussion. In the $4\times 4$ case Holtz and Sturmfels discovered that the $\mathfrak{S_{4}}\ltimes \prod_{i=1..4}\SL(V_{i})$-orbit of $h$ generates the ideal of relations, and conjectured that a similar orbit of equations, which they called the hyperdeterminantal module, generates the ideal of relations among principal minors of a symmetric matrix.
This conjecture was solved (set-theoretically) by the author \cite{Oeding_principal, oeding_thesis}. We note that in the framework of this paper, imposing that the matrix $B$ be symmetric naturally imposes a restriction from the Grassmannian $\Gr(m,2m)$ to the Lagrangian Grassmannian $\textrm{LGr}(m,2m)$, (see \cite{oeding_thesis}).

In the case that the matrix $B$ is assumed to be skew-symmetric, one considers the principal Pfaffians. It turns out that the relations among principal Pfaffians are precisely the equations of the orthogonal Grassmannian, which are analogous to the Pl\"ucker relations and were known classically, (see \cite{LM02} for a modern treatment in invariant language).

\subsection{Higher dimensional images}
It would be interesting to consider, for instance, a camera with fixed position continuously viewing a scene in time in our paradigm as $W$ a $5$-dimensional vector space and $V_{i}$ each as $4$-dimensional vector spaces (3 space dimensions and one time).  This seems to be a promising approach to understanding this and a wide variety of generalizations.

\section{Contractions and homography tensors}\label{sec:Contract}

According to Shashua and Wolf \cite{ShashuaWolf} a contraction of a quadrifocal tensor is a homography tensor \cite{Shashua00homographytensors} in the other 3 views. Moreover, this property defines the set of homography tensors. 
Similarly, the contraction of a homography tensor is an LLC mapping (a rank 2 matrix relating the epipoles in two views). We will give an invariant description of this process and show how to get 600 independent internal constraints on quadrifocal tensors.

Consider $V_{1}\o V_{2}\o V_{3}\o V_{4}$  and suppose we choose a basis of $V_{4}$, $\{x,y,z\}$. Then we may write $Q\in V_{1}\o V_{2}\o V_{3}\o V_{4}$ as 
\[Q(x,y,z) =xQ_{1} + y Q_{2} + z Q_{3},\]
 where $Q_{i}$ are the standard $3 \times 3 \times 3$ slices of the $3 \times 3 \times 3 \times 3$ tensor. By abuse of notation, we also think of $Q(x,y,z)$ as a function from $\CC^{3}$ to $V_{1}\o V_{2}\o V_{3}$, and let $x,y,z$ act as variables.

The following results come from \cite{ShashuaWolf}  and \cite{Shashua00homographytensors}.
\begin{theorem}[{\cite[Theorem~2]{ShashuaWolf}}] If $Q$ is a quadrifocal tensor then $Q(x,y,z)$ is a homography tensor for every value of $x,y,z$. 
\end{theorem}
For our purposes, the above result can be used to define the notion of ``homography tensor.''
Now choose a basis $a,b,c$ of $V_{3}$. If $H \in V_{1}\o V_{2}\o V_{3}$, then we may write $H(a,b,c) = aH_{1}+bH_{2} + cH_{3}$.  
In slightly different language  \cite[Prop.~7.2]{AholtOeding} we showed that the Zariski closure of the homography tensors is an irreducible variety (we called it $\textrm{P-Rank}^{2,2,2}$), defined (at least set-theoretically) by 30 cubic equations.

\begin{theorem}[{\cite[Theorem~1]{ShashuaWolf}}]
If $H$ is a homography tensor then $H(a,b,c)$ is an  LLC (Linear Line Complex) mapping for every value of $a,b,c$. 
\end{theorem}
For our purposes, an LLC is a skew-symmetric $3 \times 3$ matrix, which necessarily has (even) rank $\leq 2$.
Thus for all values of $a,b,c$ the matrix $H(a,b,c)$ must satisfy the constraint
$\det(H(a,b,c)) \equiv 0$. In particular, every coefficient in the cubic polynomial in $a,b,c$ must vanish.  This condition gives a basis of the ten-dimensional space of cubics. Note, the coefficient of $a^{3}$ is the determinant of the first slice, and the coefficients on $b^{3}$ and $c^{3}$ are respectively the determinants of the second and third slices.
This space also has the interpretation of $\bw{3}V_{1}^{*}\o \bw{3}V_{2}^{*}\o S^{3}V_{3}^{*}$ as a $G$-module.

We can apply the same method to quadrifocal tensors. Namely if $Q$ is a quadrifocal tensor, $Q(x,y,z)$ must satisfy all of the internal trifocal constraints for all values of $x,y,z$. In particular, $Q(x,y,z)(a,b,c)$ must be an epipolar matrix (whose entries are bi-homogeneous quadrics), and we must have the bi-homogeneous sextic polynomial (of bi-degree (3,3)) $\det(Q(x,y,z)(a,b,c)) \equiv 0$ for all values of $x,y,z,a,b,c$.  

The space of bi-degree (3,3) sextics is $100$-dimensional, and the coefficients of the expression $\det(Q(x,y,z)(a,b,c))$ provide a basis of a 100-dimensional space of cubics that vanish on the quadrifocal variety.  Note, the coefficient of $x^{3}a^{3}$ is the determinant of the first slice, and there are 8 other monomials that are the product of two cubes, the coefficients of which correspond to the determinants of the 8 other slices.  This space also has the interpretation as the $\GL(3)^{\times 4}$-module $\bw{3}V_{1}\o \bw{3}V_{2}\o S^{3}V_{3} \o S^{3}V_{4}$.  

If we interchange the roles of $V_{1},V_{2},V_{3},V_{4}$, and apply the same construction we obtain 6 non-isomorphic modules of the same format 
$\bw{3}V_{i}\o \bw{3}V_{j}\o S^{3}V_{k} \o S^{3}V_{l}$.  In particular, we find a space of 600 cubic polynomials in the ideal of the quadrifocal variety, 54 of which are determinants of $3 \times 3$ slices of a $3 \times 3 \times 3 \times 3$ tensor.  

Let $G:=\mathfrak{S}_{4}\ltimes \GL(3)^{\times 4}$. Since the quadrifocal variety is $G$-invariant, we can describe its ideal as a $G$ module. For convenience, when the $\mathfrak{S}_{4}$ symmetry is present, we write $S_{\pi_{1}}S_{\pi_{2}} S_{\pi_{3}} S_{\pi_{4}}$ for  the direct sum of
$S_{\pi_{1}}V_{1}^{*}\o S_{\pi_{2}}V_{2}^{*}\o S_{\pi_{3}}V_{3}^{*}\o S_{\pi_{4}}V_{4}^{*}$
and all non-redundant copies of it obtained by permuting the indices.

The above discussion implies the following:
\begin{prop}\label{prop:600}
Suppose $Q$ is a quadrifocal tensor. Then the 600 polynomials forming a basis of the $G$-module $S_{3}S_{3}S_{1,1,1}S_{1,1,1}$ vanish on $Q$.
\end{prop}

Because one of the contractions of a trifocal tensor also form a homography, we know that 10 cubic polynomials vanish on the set of trifocal tensors. In \cite{AholtOeding} we showed that this condition cuts out a subset of the $3 \times 3 \times 3$ tensors consisting of 4 irreducible algebraic varieties.  In order to distinguish the trifocal variety, more equations are needed. 

\begin{theorem}[ \cite{AholtOeding}]\label{thm:AO}
The ideal of the trifocal variety in $V_{1}\o V_{2}\o V_{3}^{*}$ is generated by 10 cubic, 81 quintic, and 1980 sextic polynomials. These are represented by the following $\GL(3)^{\times 3}$ modules:

\noindent$\begin{matrix}
\M_3 &=& \bw{3}V_{1}^{*}\o \bw{3}V_{2}^{*}\o S^{3}V_{3}, \\
\\
\M_{5}&=& S_{221}V_{1}^{*}\o S_{221}V_{2}^{*}\o S_{311}V_{3} &\oplus& S_{221}V_{1}^{*}\o S_{221}V_{2}^{*}\o S_{221}V_{3},  \\
\\
\M_{6}&= 
& S_{411}V_{1}^{*}\o S_{33}V_{2}^{*}  \o S_{222}V_{3} &\oplus 
& S_{33}V_{1}^{*}\o   S_{411}V_{2}^{*}\o  S_{222} V_{3}  &\oplus 
& S_{33}V_{1}^{*}\o   S_{222}V_{2}^{*}\o S_{411}V_{3} \\&\oplus 
& S_{222}V_{1}^{*}\o S_{33}V_{2}^{*} \o S_{411}V_{3}  &\oplus 
& S_{33}V_{1}^{*}\o  S_{33}V_{2}^{*} \o S_{222} V_{3} &\oplus 
& S_{33}V_{1}^{*}\o  S_{222}V_{2}^{*}\o  S_{33}V_{3} \\&\oplus 
& S_{222}V_{1}^{*}\o S_{33}V_{2}^{*} \o S_{33} V_{3} &\oplus 
& S_{33}V_{1}^{*}\o   S_{321}V_{2}^{*} \o S_{321}V_{3} &\oplus 
& S_{321}V_{1}^{*}\o  S_{33}V_{2}^{*}\o  S_{321}V_{3} .
\end{matrix}$

\end{theorem}

In the next section we work to obtain a similar statement for quadrifocal tensors.
\section{Computational results for the quadrifocal ideal}\label{sec:Rep}
The goal of this section is to give a description of the lowest degree part of the vanishing ideal for the quadrifocal variety in terms of $G$-modules.
First recall that the polynomial ring $\bigoplus_{d} S^{d}(\VVs)$ has a graded isotopic decomposition (see \cite{Landsberg-Manivel04}, for instance) with respect to $G = \mathfrak{S}_{4}\ltimes \GL(3)^{\times 4}$:
\[
S^{d}(\VVs) = \bigoplus_{\pi}S_{\pi}V^{*} \o \CC^{m_{\pi}} 
,\]
with Schur modules $S_{\pi}V^{*}:=S_{\pi_{1}}V_{1}^{*} \o S_{\pi_{2}}V_{2}^{*} \o S_{\pi_{3}}V_{3}^{*} \o S_{\pi_{4}}V_{4}^{*}$ and multiplicity space $ \CC^{m_{\pi}} $. The multiplicity space $\CC^{m_{\pi}}$ has a basis given by linear combinations of fillings of shape $\pi$.
We write $(S_{\pi_{1}}S_{\pi_{2}}S_{\pi_{3}}S_{\pi_{4}})\o \CC^{m}$ to indicate the $G$ module (occurring with multiplicity $m$) gotten by summing over all non-redundant permutations of the partitions indexing $S_{\pi}V^{*}$.

\begin{computation} Let $I_{d}$ denote the degree $d$ piece of the ideal of the quadrifocal variety, and let $G = \mathfrak{S}_{4} \ltimes \GL(3)^{\times 4}$. The symmetry assisted computations of $I_{d}$ up to degree $d\leq 9$ using random points of the quadrifocal variety are reported in Table~\ref{toDeg9}.
For $5\leq d \leq 9$ the dimensions of $I_{d}$  and the list necessary modules of minimal generators in Table~\ref{toDeg9} are correct with high probability.
\normalfont
\begin{table}
\begin{tabular}[t]{|c|r|l|r|}\hline
\begin{tabular}{l}graded\\  piece \end{tabular} &\begin{tabular}{l} $\dim I_{d}$ \end{tabular} &\begin{tabular}{l}necessary $G$ modules \\of minimal generators\end{tabular} & \begin{tabular}{l}dimension of \\ necessary mingens\end{tabular} \\
\hline
$I_{2}$ &0 &$\begin{array}{rl}\M_{2}=&0\end{array}$ &0 \\
$I_{3}$ &600 &$\begin{array}{rl}\M_{3}= &S_{3}S_{3}S_{1,1,1}S_{1,1,1}\end{array}$ &600 \\
$I_{4}$ &48,600& $\begin{array}{rl}\M_{4}=&0\end{array}$ & 0 \\
$I_{5}$ & 1,993,977 & $\begin{array}[t]{rl}\M_{5}=& S_{3,1,1}S_{3,1,1}S_{3,1,1}S_{3,1,1}\\ &\oplus   S_{2,2,1}S_{2,2,1}S_{2,2,1}S_{2,2,1}\end{array}$ & 1,377 \\
$I_{6}$ & 54,890,407 & 
$\begin{array}[t]{rl}
%(s_6*t_(3,3)*u_(3,3)+((2*s_(4,1,1)+2*s_(3,3))*t_(3,3)+s_(3,2,1)*t_(3,2,1)+2*s_(2,2,2)*t_(2,2,2))*u_(2,2,2))*v_(2,2,2)
\M_{6}=& S_{4,1,1}S_{3,3}S_{2,2,2}S_{2,2,2}\o \CC^{2}\\
&\oplus S_{3,3}S_{3,3}S_{2,2,2}S_{2,2,2}\o \CC^{2}\\
&\oplus S_{3,2,1}S_{3,2,1}S_{2,2,2}S_{2,2,2}\\
&\oplus S_{2,2,2}S_{2,2,2}S_{2,2,2}S_{2,2,2}\o \CC^{2}\\
&\oplus S_{6}S_{3,3}S_{3,3}S_{2,2,2}
\end{array}$
&37,586 \\
$I_{7}$ &1,140,730,128 &$ \begin{array}{rl}\M_{7}=&  0 \end{array}$& 0 \\
$I_{8}$ &
18,940,147,947
%??18,051,062,139
%18,566,037,864 -- old
 & $\begin{array}[t]{rl}\M_{8}= & S_{4,4}S_{4,4}S_{4,4}S_{4,2,2} \otimes \CC^{2}\end{array}$ &
162,000
\\
$I_{9}$ & $\geq$ 223,072,284,455
%old:?188,850,321,637 
& $ \begin{array}[t]{rl}\M_{9}\geq &  S_{5, 4} S_{5, 4} S_{5, 4}S_{4, 3, 2} \\ & \oplus S_{5, 4} S_{5, 4} S_{5, 4} S_{5, 2, 2}\end{array}$ &$ \geq$ 3,087,000\\
\hline
\end{tabular}
\caption{The dimension and isotypic description of the ideal of the quadrifocal variety up to degree 9.}\label{toDeg9}
\end{table}

\end{computation}
\emph{Description of computation:} 
The qualifier ``with high probability'' refers to the fact that we computed the ideal on random subsets of points from the quadrifocal variety, so in principle we could have chosen a set of points in special position and would have over counted the dimension of the ideal. The set of points in such special position, however, being of lower dimension, has measure zero in the quadrifocal variety, so we say the results hold with high probability.  The proof techniques we use for the degree $\leq 4$ computations, however are valid unconditionally.

The content of this computation is two-fold.  First we computed (in Maple) the entire ideal degree by degree, up to degree 8 and partially in degree 9, making use of the isotypic decomposition of the polynomial ring.  Second, we checked for representation-theoretic certificates for necessity of minimal generators using \texttt{SchurRings} in Macaulay2.
In the ancillary files associated with the arXiv version of this manuscript we provide a minimal set of data necessary to check our work. This includes basic maple scripts, fillings that yield a basis of each isotypic component, and the resulting modules in the ideal. 
Here is a summary of these computations.

Every irreducible $G$-module has the property that it is the vector space spanned by the $G$-orbit of a so-called \emph{highest weight vector} (see \cite{LandsbergTensorBook}).  We used Young symmetrizers to obtain a basis of each multiplicity space $\CC^{m_{\pi}}$ in the isotypic decomposition of the polynomial ring above. Then we computed the subspace of $\CC^{m_{\pi}}$ that vanished on the quadrifocal tensors. This algorithm was used and outlined in \cite{AholtOeding, BatesOeding} for example.  Since we learned this algorithm from the paper of Landsberg and Manivel \cite{Landsberg-Manivel04}, we call this algorithm the \emph{Landsberg-Manivel-algorithm} or the \emph{LM-algorithm} for short.

The output of this symmetry-enhanced polynomial interpolation computation is a $G$-module description of the ideal in each degree.
The dimensions of these modules give the beginning of the Hilbert function of the quadrifocal ideal, and are also reported in Table~\ref{toDeg9}.
We will use notation of \texttt{SchurRings}:
The ring $Sym^{\bullet}( \VVs) $ is regarded as a tower of rings. The variables in each ring are represented as $s_{\pi}$, (respectively $t_{\pi}$, $u_{\pi}$, $v_{\pi}$), for partitions $\pi$. The correspondence between the two notations is
\[{m} s_{\pi_{1}}t_{\pi_{2}}u_{\pi_{3}}v_{\pi_{4}} \leftrightarrow S_{\pi_{1}}V_{1}^{*} \o S_{\pi_{2}}V_{2}^{*} \o S_{\pi_{3}}V_{3}^{*} \o S_{\pi_{4}}V_{4}^{*} \o \CC^{m}.\]

In degree $3$ our application of the LM-algorithm found the following modules:
{\small \[I_{3}=
({s}_{(1,1,1)} {t}_{(1,1,1)} {u}_{3}+\blue({s}_{(1,1,1)} {t}_{3}+{s}_{3} {t}_{(1,1,1)}\blue) {u}_{(1,1,1)}) {v}_{3}+(\blue({s}_{(1,1,1)} {t}_{3}+{s}_{3} {t}_{(1,1,1)}\blue) {u}_{3}+{s}_{3} {t}_{3}
      {u}_{(1,1,1)}) {v}_{(1,1,1)}.\]}
To save space, we only record those modules up to the $\mathfrak{S}_{4}$-action:
 \[\M_{3}:={s}_{3} {t}_{3} {u}_{(1,1,1)} {v}_{(1,1,1)}\]
This module corresponds to the same $600$ polynomials in Proposition~\ref{prop:600}. Non-vanishing of polynomials on random points of a variety implies non-vanishing for the entire variety, so this $600$ dimensional vector space of cubics are the only cubics vanishing on the quadrifocal variety. Therefore, the $d\leq 3$ computations hold with no ``high probability'' qualifier.

The modules we found (by the LM-algorithm) in degree 4 are (up to the $\mathfrak{S}_{4}$-action):
\[I_{4}=({s}_{4} {t}_{4}+\blue({s}_{4}+{s}_{(3,1)}\blue) {t}_{(3,1)}) {u}_{(2,1,1)} {v}_{(2,1,1)}\]

Because $81*600 = 48600$, we guess that all the equations in $I_{4}$ come from linear combinations of products of linear forms with  the $600$ cubics in $I_{3}$ and thus we would guess that there are no new generators in degree 4. We would like to prove this using representation theory. 

The multiplication in the ring $S^{\bullet} (\VVs)$ is just the usual polynomial multiplication. It is not easy to determine the isotopic version of multiplication.
However, we can get a lower bound on the modules of minimal generators in our ideal. The basic idea is the following.
The multiplication in the ring $S^{\bullet}(\VVs)$ is also the restriction of the multiplication in the tensor ring $(\VVs)^{\otimes}$. The tensor product of two representations of the form $S_{\pi}V^{*}$ is obtained by an iteration of the Littlewood-Richardson rule.  

Using the package \texttt{SchurRings} \cite{SchurRings} we computed the tensor product
\[
I_{3}\o( \VVs)
\]
and found that it exactly coincides with $I_{4}$. This is an indication that there are probably no new generators in degree $4$, however it could be that the modules resulting in $I_{3}\o \VVs$ are not actually in $S^{4}(\VVs)$. 
The only way that there could be new minimal generators in degree 4 is if there were already some syzygies amongst the degree 3 equations. 
At least in degree 4 it is possible to look at the highest weight vectors and check that they are in the ideal generated by the 600 cubics in Macaulay2.  On the other hand, we can argue in a less computationally intensive way by using the following special case of a more general idea from \cite{RaicuProducts}, which was employed in \cite{RaicuGSS, OedingRaicu}.
\begin{lemma}\label{removeBox}
Suppose $F_{\pi}:=(F_{\pi_{1}},F_{\pi_{2}},F_{\pi_{3}},F_{\pi_{4}})$ is a filling using the ordered alphabets $(\A_{1},\A_{2},\A_{3},\A_{4})$ giving a nonzero realization of the module $S_{\pi}V^{*}$ in $S^{d}(\VVs)$. If $F_{\mu}:=(F_{\mu_{1}},F_{\mu_{2}},F_{\mu_{3}},F_{\mu_{4}})$ is a filling giving a nonzero realization of the module  $S_{\mu}V^{*}$ in $S^{d-1}(\VVs)$ and $F_{\mu}$ may be obtained from $F_{\pi}$ by respectively deleting the last used letter in each of the alphabets $(\A_{1},\A_{2},\A_{3},\A_{4})$, then $S_{\pi}V^{*}$ is in the ideal generated by $S_{\mu}V^{*}$.
\end{lemma}
%\begin{proof}
%Compute a highest weight vector of $S_{\mu}$ by applying the Young symmetrizer algorithm to the specified fillings, then use the Littlewood-Richardson rule (actually the Pieri rule) to compute $S_{\mu}\o S_{1}$. 
%Note that $S_{\pi}$ must occur in the tensor product by assumption.  The highest weight vector in $S_{\mu} \otimes S_{1}$ of weight $\pi$ must occur in $S_{\pi}$ and is non-zero.  
%\end{proof}

Now we apply Lemma~\ref{removeBox} to see that $S_{211}S_{211}S_{31}S_{31}$ is in the ideal generated by $S_{111}S_{111}S_{3}S_{3}$. In this case everything occurs with multiplicity one, so our work is much easier.
The filling
\[
\young(14,2,3)\o \young(14,2,3)\o \young(123,4)\o \young(123,4)
\]
produces a realization of a copy of $S_{31}S_{31}S_{211}S_{211}$ in the ideal of the quadrifocal variety. Notice that by removing the fourth letter we obtain the filling
\[
\young(1,2,3)\o \young(1,2,3)\o \young(123)\o \young(123)
,\]
which produces a nonzero copy of $S_{111}S_{111}S_{3}S_{3}$.

The same argument may be applied to  $S_{4}S_{4}S_{211}S_{211}$, with the same result.
Therefore, none of the modules in degree 4 are minimal generators of the ideal of the quadrifocal variety.

In degree 5 we note that $1929501-600*3321 = 1,377$. This means that the degree 3 equations cannot generate all of the ideal in degree 5, and there must be at least $1,377$ new minimal generators in degree 5. Moreover, if there are no degree 2 syzygies amongst the degree 3 equations then the space of minimal generators in degree 5 would be precisely $1,377$-dimensional.  In principle one could try to use a degree-limited Gr\"obner basis computation in  Macaulay2  to check if each highest weight vector of each of the modules in $I_{5}$ is actually in $\langle I_{3}\rangle_{5}$, but memory limitations become a problem.

However, we can argue by comparing multiplicities of isotypic components. 
In degree 5 the LM-algorithm produced the following modules (up to the $\mathfrak{S}_{4}$-action):
{\small \begin{multline*}I_{5} =  \red({s}_{5} {t}_{5}+\blue({s}_{5}+2 {s}_{(4,1)}\blue) {t}_{(4,1)}+\blue({s}_{5}+{s}_{(4,1)}+{s}_{(3,2)}\blue) {t}_{(3,2)}+\blue(3 {s}_{(4,1)}+7 {s}_{(3,1,1)}\blue) {t}_{(3,1,1)}\red) {u}_{(3,1,1)} {v}_{(3,1,1)}\\
+\green(\red(\blue({s}_{5}+2{s}_{(4,1)}\blue) {t}_{(4,1)}+{s}_{(4,1)} {t}_{(3,2)}+\blue(2 {s}_{5}+2 {s}_{(4,1)}+2 {s}_{(3,2)}\blue) {t}_{(3,1,1)}\red) {u}_{(3,1,1)}\\
+\red({s}_{5} {t}_{5}+\blue({s}_{5}+2 {s}_{(4,1)}\blue){t}_{(4,1)}+\blue({s}_{5}+{s}_{(4,1)}+{s}_{(3,2)}\blue) {t}_{(3,2)}
+\blue({s}_{(4,1)}+{s}_{(3,1,1)}\blue) {t}_{(3,1,1)}
+{s}_{(2,2,1)} {t}_{(2,2,1)}\red)
 {u}_{(2,2,1)}\green) {v}_{(2,2,1)}.
\end{multline*}
}

Using \texttt{SchurRings} we found that the only modules left in the difference between $I_{5}$ and $I_{4}\otimes \VV$ are 
 \[\M_{5}:={s}_{(3,1,1)} {t}_{(3,1,1)} {u}_{(3,1,1)} {v}_{(3,1,1)}+{s}_{(2,2,1)} {t}_{(2,2,1)} {u}_{(2,2,1)} {v}_{(2,2,1)}.\] From this computation we know that these two summands must be among the minimal generators.  We do not know if there are any other minimal generators in degree 5. 
 \begin{remark}While this \texttt{ShurRings} computation ignores the structure of the multiplicity spaces in the ideals,  to first approximation it tells if there is any representation-theoretic reason for modules to be among the minimal generators.
\end{remark}
In degree 6 the LM-algorithm produced the following modules (up to the $\mathfrak{S}_{4}$-action) $I_{6}=$
{\small \begin{multline*}
\red{(}s_{6}t_{6}+\blue{(}s_{6}+2s_{(5,1)}\blue{)}t_{(5,1)}+\blue{(}s_{6}+2s_{(5,1)}+2s_{(4,2)}\blue{)}t_{(4,2)}+\blue{(}3s_{(5,1)}+3s_{(4,2)}+10s_{(4,1,1)}\blue{)}t_{(4,1,1)}\red{)}u_{(4,1,1)}v_{(4,1,1)}\\
+\red(\blue{(}s_{6}+s_{(5,1)}+s_{(4,2)}\blue{)}t_{(4,1,1)}u_{(4,1,1)}+s_{(4,1,1)}t_{(4,1,1)}u_{(3,3)}\red)v_{(3,3)}\\
+\green(\red(\blue{(}s_{6}+3s_{(5,1)}\blue{)}t_{(5,1)}+\blue{(}s_{6}+3s_{(5,1)}+3s_{(4,2)}\blue{)}t_{(4,2)}+\blue(2s_{6}+7s_{(5,1)}+7s_{(4,2)}+12s_{(4,1,1)}\blue{)}t_{(4,1,1)}\red)u_{(4,1,1)}\\
+\blue{(}s_{(5,1)}+s_{(4,2)}+2s_{(4,1,1)}\blue{)}t_{(4,1,1)}u_{(3,3)}\\
+\red(s_{6}t_{6}+\blue(2s_{6}+6s_{(5,1)}\blue{)}t_{(5,1)}+\blue(2s_{6}+6s_{(5,1)}+6s_{(4,2)}\blue{)}t_{(4,2)}+\blue(3s_{6}+10s_{(5,1)}+10s_{(4,2)}+16s_{(4,1,1)}\blue{)}t_{(4,1,1)}\\
+\blue{(}s_{6}+2s_{(5,1)}+2s_{(4,2)}+4s_{(4,1,1)}+s_{(3,3)}\blue{)}t_{(3,3)}\\
+\blue(3s_{6}+12s_{(5,1)}+12s_{(4,2)}+16s_{(4,1,1)}+3s_{(3,3)}+24s_{(3,2,1)}\blue{)}t_{(3,2,1)}\red)u_{(3,2,1)}\green)v_{(3,2,1)}\\
+\green(\red{(}s_{(5,1)}t_{(5,1)}+\blue{(}s_{(5,1)}+s_{(4,2)}\blue{)}t_{(4,2)}
+\blue{(}s_{6}+3s_{(5,1)}+3s_{(4,2)}+6s_{(4,1,1)}\blue{)}t_{(4,1,1)}\red)u_{(4,1,1)}\\
+\red{(}s_{(4,1,1)}t_{(4,1,1)}+s_{6}t_{(3,3)}\red)u_{(3,3)}+\red(\blue{(}s_{6}+3s_{(5,1)}\blue{)}t_{(5,1)}+\blue{(}s_{6}+3s_{(5,1)}+3s_{(4,2)}\blue{)}t_{(4,2)}\\
+\blue{(}s_{6}+4s_{(5,1)}+4s_{(4,2)}+5s_{(4,1,1)}\blue{)}t_{(4,1,1)}+\blue{(}s_{(5,1)}+s_{(4,2)}+s_{(4,1,1)}\blue{)}t_{(3,3)}
\\
+\blue{(}s_{6}+4s_{(5,1)}+4s_{(4,2)}+6s_{(4,1,1)}+3s_{(3,3)}+8s_{(3,2,1)}\blue{)}t_{(3,2,1)}\red)u_{(3,2,1)}\\
+\red{(}s_{6}t_{6}+2s_{(5,1)}t_{(5,1)}+\blue{(}s_{6}+s_{(5,1)}+2s_{(4,2)}\blue{)}t_{(4,2)}\\
+\blue{(}2s_{(5,1)}+2s_{(4,2)}+2s_{(4,1,1)}\blue{)}t_{(4,1,1)}+\blue{(}s_{(5,1)}+3s_{(4,1,1)}+3s_{(3,3)}\blue{)}t_{(3,3)}\\
+\blue{(}s_{(5,1)}+s_{(4,2)}+s_{(4,1,1)}+4s_{(3,2,1)}\blue{)}t_{(3,2,1)}+3s_{(2,2,2)}t_{(2,2,2)}\red)u_{(2,2,2)}\green)v_{(2,2,2)}
.      \end{multline*}
      }
Applying the same \texttt{SchurRings} test we find the following necessary minimal generators:
{\small \begin{multline*}\M_{6}:=
\green({s}_{6} {t}_{(3,3)} {u}_{(3,3)}+\red(\blue(2 {s}_{(4,1,1)}+2 {s}_{(3,3)}\blue) {t}_{(3,3)}+{s}_{(3,2,1)} {t}_{(3,2,1)}+2 {s}_{(2,2,2)} {t}_{(2,2,2)}\red) {u}_{(2,2,2)}\green) {v}_{(2,2,2)}
      .\end{multline*}
}
%\red{The new computation agrees with the previous computation:
%\begin{multline*}
%2v_{2,2,2}u_{2,2,2}t_{2,2,2}s_{2,2,2}
%+v_{2,2,2}u_{3,3}s_6t_{3,3}
%+2v_{2,2,2}u_{2,2,2}s_{4,1,1}t_{3,3} \\
%+2v_{2,2,2}u_{2,2,2}t_{3,3}s_{3,3}
%+v_{2,2,2}u_{2,2,2}t_{3,2,1}s_{3,2,1}
%\end{multline*}
%}

%This module together with all the copies added by the $\mathfrak{S}_{4}$-symmetry, has dimension 37,586.
% and the Schur functor notation is given in the statement of the theorem.
%%% the previous case when I didn't account for combinations of M5 and I3
%\M_{6}:=4 {s}_{(3,2,1)} {t}_{(3,2,1)} {u}_{(3,2,1)} {v}_{(3,2,1)}\\
%+(3 {s}_{(4,1,1)} {t}_{(4,1,1)} {u}_{(4,1,1)}+{s}_{6} {t}_{(3,3)} {u}_{(3,3)}
%+({s}_{(4,1,1)}+{s}_{(3,3)}+4 {s}_{(3,2,1)})
%      {t}_{(3,2,1)} {u}_{(3,2,1)}\\
%      +(({s}_{(5,1)}+{s}_{(4,2)}) {t}_{(4,1,1)}+(3 {s}_{(4,1,1)}+2 {s}_{(3,3)}) {t}_{(3,3)}+2 {s}_{(3,2,1)} {t}_{(3,2,1)}+2 {s}_{(2,2,2)} {t}_{(2,2,2)}) {u}_{(2,2,2)})
%      {v}_{(2,2,2)}.
These modules could not come from any of the lower degree parts of the ideal, so they are among the minimal generators in degree 6, however there could be more minimal generators that we haven't accounted for if there were syzygies among the lower degree generators.

The space of degree 7 polynomials $S^{7}(\VVs)$ decomposes as a sum of 288 isotypic components (up to $\mathfrak{S}_{4}$ symmetry) with dimension of multiplicity spaces as large as 301. Thus applying the LM-algorithm directly in degree 7 was more challenging and we were forced to exploit a parallelism, following a process similar to that in \cite{OedingSam}.   We ran the LM-algorithm in a separate instance of Maple for each of the 288 isotypic components. This allowed us to distribute the computation. Running  in parallel the degree 7 computation took approximately 12 hours on our two servers with respectively 40 and 24 processors. We found 201 modules occurring non-trivially the ideal. The   \texttt{SchurRings}  test yielded no new necessary minimal generators.

%\red{Update: After fixing the number of points we evaluate at (decreasing the multiplicities of the modules in lower degrees as well as this degree), I still find no new mingens in degree 7 implied by this calculation}

The parallelized LM-algorithm for degree 8 involved 619 Schur modules with multiplicity spaces up to 608-dimensional. 
Parts of this computation took approximately 20 days on our servers. We found 453 modules occurring in $I_{8}$ with nonzero multiplicity.
The \texttt{ShurRings} test only produced one module occurring in $I_{8}$ with greater multiplicity than can occur in $I_{7}\o \VVs$, namely
\[
\M_{8}:=
2s_{(4,4)}t_{(4,4)}u_{(4,4)}v_{(4,2,2)}
.\]
That this module is not in the ideal generated by $\M_{3}+\M_{5}+\M_{6}$ can be seen for shape reasons alone.
The Schur module $\M_{8}$ is indexed by a multi-partition containing three partitions of shape $(4,4)$, but no module the list of known minimal generators is indexed by quadruple of partitions with a subset of 3 of them that each fit in a box of dimensions $2\times 4$.  Therefore $\M_{8}$ must occur with multiplicity 2 in the minimal generators in the ideal.
% $\M_{8}$ is the only module that occurred with greater multiplicity in $I_{8}$ than in $I_{7}\o \VVs$. 

In degree 9 there are 1205 isotypic components, with multiplicity spaces as high as 2226-dimensional.  After approximately one month of computational time on our servers we were able to evaluate 1158 of these isotypic components on the quadrifocal variety, 951 of which occurred non-trivially in the ideal. The computations that finished were those modules for which the multiplicity was low (under around 300). 

\begin{remark}The main obstruction to completing the degree 9 computation is a lack of available machines able to perform Maple computations.  More specifically, if a given isotypic component has multiplicity $m$ in the polynomial ring, we must populate an $m\times m$ matrix with $m^{2}$ evaluations (from the Young symmetrizer algorithm), once to verify we have found a basis of the space, and a second time to evaluate this basis on points of the quadrifocal variety. Each of these evaluations could be done in a separate instance. So in the case of multiplicity 2226, we must perform approximately 10 million evaluations, which could potentially be done on a separate processor if were to take fuller advantage of the parallel nature of this problem.  One of the largest multiplicity spaces in degree 9 we were able to handle was 456 dimensional, and the computation in multiplicity 2226 case is approximately 23 times larger than this.
Since this scale of resources is not currently available to us, it seems unlikely that we will be able to complete the degree 9 computation. 
\end{remark}

Though we were only able to partially compute $I_{9}$, we were able to obtain two new modules in degree 9 that must be among the minimal generators of the quadrifocal ideal:
\[
\widetilde \M_{9}:=s_{(5, 4)} t_{(5, 4)} u_{(5, 4)} v_{(4, 3, 2)} +
s_{(5, 4)} t_{(5, 4)} u_{(5, 4)} v_{(5, 2, 2)}
.\]
%which has dimension  $3,087,000$.
These two modules occur with multiplicity 3 in $I_{9}$, but instances of these modules coming from $I_{8}\o\VVs$ occurred with multiplicity 2. Therefore, these modules must occur with multiplicity at least 1 among the minimal generators of $I$.

The dimension counts in Table~\ref{toDeg9} are a straightforward application of a hook length formula, which is also implemented in \texttt{SchurRings} for example.% $\diamondsuit$

%\end{proof}
% had to be careful to take the proper symmetrization of each...
% this is still odd... double check...
We wondered if Table~\ref{toDeg9} might be a complete list of minimal generators for the quadrifocal ideal.  
We performed further experiments with modules in degrees 9, 10, and 11 having low multiplicity. 
We found that  $S_{5, 5} S_{5, 5} S_{5, 5}S_{4, 3, 3}$ occurred with multiplicity one greater in $I_{10}$ than what could come from the part of $I_{9}$ that we were able to compute.  This was the only new module we found in low multiplicity in degree 10.  Further experiments suggest that there may be many more new modules of minimal generators with low multiplicity in degree 11. At present we are not confident enough to form a conjecture as to what happens in higher degree, nor in which degree the last set of minimal generators must occur. 

%\caption{The low multiplicity part of the ideal of the quadrifocal variety up to degree 14.}\label{toDeg9}
%\end{table}

%\begin{conjecture}\label{conj}
%The ideal of the quadrifocal variety is minimally generated by 
%$600$ cubics, $1,377$ quintics, $37,586$ sextics, $162,000$ octics and $3,087,000$ nonics.
%\end{conjecture}
\begin{remark}
In a previous draft we had incorrectly computed the number of equations in the ideal because of an unfortunate programming error.
\end{remark}

\section*{Acknowledgements}
The author would like to thank Bernd Sturmfels for suggesting this problem and Claudiu Raicu for help with \texttt{SchurRings}. F\'ed\'eric Holweck also provided useful comments. 
The author is also grateful to the S. Korean National Institute for Mathematical Sciences (NIMS) and the Simons Institute for the Theory of Computing for their generous support while this work was carried out.  The author also thanks an anonymous referee, whose suggestions improved the exposition of this manuscript.

Part of the computational work reported in this work was performed on the Auburn CASIC High Performance Computing Cluster. The author is grateful for the availability and access to this resource.

\bibliographystyle{amsplain}
\bibliography{/Users/oeding/Dropbox/BibTeX_bib_files/main_bibfile}

\end{document}